\input amstex
\input amsppt.sty
\magnification\magstep1

\def\as{\operatorname{asdim}}
\def\G{\Gamma}

\def\Z{\text{\bf Z}}

\def\Ucal{\Cal{U}}
\def\Vcal{\Cal{V}}

\hoffset= 0.0in
\voffset= 0.0in
\hsize=32pc
\vsize=38pc
\baselineskip=24pt
\NoBlackBoxes
\topmatter
\author
J. Smith
\endauthor

\title
The Asymptotic Dimension of the First Grigorchuk Group is Infinity
\endtitle

\thanks 
\endthanks

\address University of Florida, Department of Mathematics, P.O.~Box~118105,
358 Little Hall, Gainesville, FL 32611-8105, USA
\endaddress


\email  justins\@math.ufl.edu
\endemail

\endtopmatter

\document

 The first Grigorchuk group is described in 
\cite{Gr1}, \cite{Gr2}, and 
\cite{Ha}. This Grigorchuk group, which we will denote by $\G$, 
has many interesting properties.  It is a finitely generated $2-$group with intermediate 
growth, whose word problem is solvable, and 
which does not admit a finite dimensional linear representation that is faithful.   
Also, $\G$ and $\G \times \G$ are {\it commensurable}, which means that $\G$ and $\G \times \G$ 
have subgroups of finite index which are isomorphic. 
A detailed exposition can be found in \cite{Ha}. 

We prove that $\G$ has asymptotic dimension infinity, $\as \G=\infty$. 
If one excludes Gromov's ``random groups'' \cite{Gro1}, all previously
known examples of groups $G$ with $\as G =\infty$ 
are based on the fact that $G$
has a free abelian subgroup of arbitrary large rank.
The Grigorchuk group is of different nature:
since $\G$ is a $2-$group, it does not have a (nontrivial) 
free abelian subgroup.

 Let $(X,d_X)$ and $(Y,d_Y)$ be metric spaces.  We say that a map $f:X \to Y$  
is proper if the preimage of each bounded set is bounded; it is bornologous if, for all 
$R > 0$, there is an $S>0$ such that $d_Y (f(x),f(y)) < S$ whenever 
$d_X (x,y) < R$; a map is called coarse if it is both proper and 
bornologous.  If $S$ is a set, then 
$f,g:S \to X$ are said to be close if $ \sup_{s \in S} d(f(s),g(s)) < \infty$.  
A coarse map $f: X \to Y$ is said to be a coarse equivalence if 
there is a coarse map $g: Y \to X$ such that $g\circ f$ is close to 
$\text{id}_X$ and $f \circ g$ is close to $\text{id}_Y$.  A coarse map $f:X \to Y$ 
is said to be a coarse embedding if $f:(X,d_X) \to (f(X), d_Y|_{f(X)})$ is a coarse equivalence.   
In particular, an isometric embedding is a coarse embedding.  Also, 
if $f_i : (X_i, d_{X_i}) \to (Y_i, d_{Y_i})$ ($i=1,2$) are coarse equivalences, then 
so is 
$f_1 \times f_2 : (X_1 \times X_2, \delta_1) \to  (Y_1 \times Y_2, \delta_2)$, where 
$\delta_1$ and $\delta_2$ are the corresponding sum metrics.    

Restricting attention to finitely generated groups, since any two word 
metrics (say $d_1$ and $d_2$) on a finitely generated group 
$G$ are equivalent, the map $\text{id}_G:(G,d_1) \to (G,d_2)$ is a coarse 
equivalence.  When 
$G$ and $H$ are groups equipped with word metrics $d_G$ and $d_H$, then the sum metric 
$d_G + d_H$ on $G \times H$ is also a word metric (for the obvious generating set).  If 
$H \leq G$ is a subgroup of finite index of (the finitely generated group) $G$, 
then the inclusion map $H \to G$ is a coarse equivalence.  
Also, an isomorphism is a coarse equivalence.  
See \cite{Roe} for further details.   

\definition{Definition \cite {Gro2}} A metric space $(X,d)$ is said 
to have $\as X \leq n$, if for each $R > 0$, there is an $S > 0$ 
and $R-$disjoint, $S-$bounded families 
$\Ucal_0, \Ucal_1, \ldots, \Ucal_n$ of subsets of $X$ such that 
$\Ucal := \cup_i \Ucal_i$ 
is a cover of $X$.  \enddefinition

We say that a family $\Vcal$ of subsets of $X$ is $R-$disjoint 
if $d(U,V) \geq R$ for all $U,V \in \Vcal$ with $U \neq V$; $\Vcal$ is said to 
be $S-$bounded if $\text{diam}V \leq S$ for all $V \in \Vcal$.  One can show that 
coarsely equivalent spaces have the same asymptotic dimension.  Thus, for finitely 
generated groups, the asymptotic dimension of the group does not depend on the choice 
of the word metric.  

\definition{Definition} Two groups $\G_1$ and $\G_2$ are {\it commensurable} if there exist 
subgroups $H_1 \leq \G_1$ and $H_2 \leq \G_2$, each of finite index, 
such that $H_1$ and $H_2$ are isomorphic.  \enddefinition 

By the comments above, $\as \G_1 = \as \G_2$ if $\G_1$ and $\G_2$ are commensurable.        

\proclaim{Theorem}  Let $G$ be a finitely generated, infinite group which is 
commensurable with its square $G\times G$.  Then 
$\as G =\infty$. \endproclaim 
\demo{Proof} We first show that $G^n$ is coarsely equivalent to 
$G$ for all $n\geq 1$.  Proceeding inductively (the $n=1$ case is immediate), we assume 
$G^n$ is coarsely equivalent to $G$.  But $G^{n+1}$ is coarsely equivalent to $G^n \times G$, which 
in turn is coarsely equivalent to $G \times G$, and so by hypothesis 
$G^{n+1}$ is equivalent to $G$.  This proves that $\as G^n = \as G$ for all $n \geq 1$.  

Also, by Exercise IV.A.12 of \cite{Ha}, there is an isometric embedding 
$f: \Z \to G$, where $G$ is taken with a word metric.  Thus, for each $n \geq 1$, 
we have an isometric embedding 
$f\times f \times \cdots \times f: \Z^n \to G^n$, where 
we take the sum metrics on $\Z^n$ and $G^n$.  Since an isometric embedding is 
a coarse embedding, we have that $\as G^n \geq \as \Z^n = n$.  
Thus, $\as G \geq n$ for all $n$.  \qed   
\enddemo

\proclaim{Corollary} Let $\Gamma$ be the Grigorchuk group.  Then $\as \Gamma =\infty$. \endproclaim 
\demo{Proof}  $\Gamma$ is finitely generated by definition.  Proposition VIII.14 and Corollary VIII.15 from 
\cite{Ha} show that $\G$ satisfies the hypotheses of the theorem. \enddemo   

It is interesting to note that $\as \G = \infty$, yet $\G$ does not contain an 
isomorphic copy of $\Z^n$.  However, $\Z^n$ does coarsely embed into $\G$. 

Finally, if one has a finitely generated group which is known to be commensurable with its 
square, then the asymptotic dimension is either 0 or infinity, depending 
on whether the group is finite or infinite.


\Refs \widestnumber\key{B-D1}


\ref\key Gr1
\by R. I. Grigorchuk
\paper Bernside's problem on periodic groups  
\jour Functional Anal. Appl. 
\yr 1980 
\vol 14 
\pages 41 - 43
\endref

\ref\key Gr2
\by R. I. Grigorchuk
\paper An example of a finitely presented amenable group not belonging to the class EG   
\jour Sbornik Math. 
\yr 1998 
\vol 189:1
\pages 75-95
\endref


\ref\key Gro1
\by M. Gromov
\paper Random walk in random groups 
\jour Geom. Funct. Anal. 
\yr  2003
\vol 13
\pages 73-146 
\issue 1
\endref


\ref\key Gro2
\by M. Gromov
\paper Asymptotic invariants of infinite groups
\inbook Geometric Group Theory
\yr 1993
\bookinfo London Math. Soc. Lecture Note Ser., 182
\eds G. Niblo and M. Roller 
\publ Cambridge University Press 
\endref


\ref\key Ha
\by Pierre de la Harpe 
\book Topics in Geometric Group Theory 
\bookinfo Chicago Lectures in Mathematics 
\yr 2000
\publ The University of Chicago Press 
\endref


\ref\key Roe
\by J. Roe
\book Lectures on coarse geometry 
\bookinfo University Lecture series 
\vol 31 
\yr 2003
\publ AMS
\endref

\endRefs
\enddocument